\documentclass[a4paper,12pt]{amsart}

\usepackage{amsthm}
\usepackage{amssymb, amsmath}
\usepackage{graphics}
\usepackage[dvips]{graphicx}

\usepackage{latexsym}
\usepackage{amscd}
\usepackage{float}
\usepackage{psfrag}
\usepackage{pst-node}
\usepackage{pst-plot}

\newtheorem{thm}{Theorem}[section]

\newtheorem{lmm}[thm]{Lemma}

\newtheorem{prp}[thm]{Proposition}

\theoremstyle{definition}
\newtheorem{dfn}[thm]{Definition}

\theoremstyle{remark}
\newtheorem*{rem}{Remark}

\newdimen\argwidth
\def\db[#1\db]{%
  \setbox0=\hbox{$#1$}\argwidth=\wd0
  \setbox0=\hbox{$\left[\box0\right]$}
    \advance\argwidth by -\wd0
  \left[\kern.3\argwidth\box0 \kern.3\argwidth\right]}

\newcommand{\Image}{\operatorname{Im}}

\newcommand{\xto}{\xrightarrow}

\newcommand{\Hom}{\mathop{\mathrm{Hom}}\nolimits}

\newcommand{\Spec}{\operatorname{Spec}}
\newcommand{\GL}{GL}
\newcommand{\SL}{SL}
\newcommand{\bCx}{\bC^{\times}}

\newcommand{\Aut}{\operatorname{Aut}}

\newcommand{\module}{\operatorname{mod}}

\newcommand{\bC}{\ensuremath{\mathbb{C}}}

\newcommand{\bP}{\ensuremath{\mathbb{P}}}
\newcommand{\bQ}{\ensuremath{\mathbb{Q}}}
\newcommand{\bR}{\ensuremath{\mathbb{R}}}

\newcommand{\bT}{\ensuremath{\mathbb{T}}}

\newcommand{\bZ}{\ensuremath{\mathbb{Z}}}

\newcommand{\scG}{\ensuremath{\mathcal{G}}}

\newcommand{\scI}{\ensuremath{\mathcal{I}}}

\newcommand{\scM}{\ensuremath{\mathcal{M}}}

\newcommand{\scO}{\ensuremath{\mathcal{O}}}

\newcommand{\nin}{\not \in}

\newcommand{\Etilde}{\widetilde{E}}
\newcommand{\Gtilde}{\widetilde{G}}
\newcommand{\ptilde}{\widetilde{p}}
\newcommand{\Gammatilde}{\widetilde{\Gamma}}
\newcommand{\scMtilde}{\widetilde{\scM}}

\newcommand{\Cone}{\operatorname{Cone}}
\newcommand{\perfmat}{\operatorname{Perf}}

\newcommand{\vtilde}{\tilde{v}}
\newcommand{\Utilde}{\widetilde{U}}
\newcommand{\bTtilde}{\widetilde{\bT}}

\title[Moduli spaces associated with dimer models]
{On moduli spaces of quiver representations
associated with dimer models}

\dedicatory{Dedicated to Professors Iku Nakamura and Eiichi Sato on their sixtieth birthdays}           

\author{Akira Ishii and Kazushi Ueda}

\begin{document}

\maketitle

\section{Introduction}

Dimer models are introduced by string theorists
to study four-dimensional $N = 1$ superconformal field theories.
See e.g. a review by Kennaway \cite{Kennaway_BT}
and references therein for a physical background.
A dimer model is a bipartite graph on a real two-torus
which encodes the information of a quiver with relations.
A typical example of such a quiver is the McKay quiver determined by
a finite abelian subgroup $G$ of $\SL(3, \bC)$ (see \cite{Reid, Ueda-Yamazaki1}).
In this case, the moduli space of representations of the McKay quiver
(for the dimension vector $(1, 1, \dots, 1)$)
coincides with the moduli space of $G$-constellations
considered in \cite{Craw-Ishii}.
For a generic choice of a stability parameter $\theta$,
the moduli space of $G$-constellations is a crepant resolution
of the quotient singularity $\bC^3/G$
and the derived category of coherent sheaves on the moduli space
is equivalent to the derived category of finitely-generated modules
over the path algebra of the McKay quiver.
It is expected that
these kinds of statements can be generalized
to the case of dimer models
that are ``consistent'' in the physics context,
which should be called ``brane tilings".
In this note,
we discuss a slightly weaker notion of
{\em non-degenerate dimer models},
which is strong enough to ensure that
the moduli space is a crepant resolution
of the three-dimensional toric singularity
determined by the Newton polygon of the characteristic polynomial
(see Theorem \ref{th:main}).
We expect that one has to impose further conditions
to prove the derived equivalence.

For the proof, we use a generalization of the description
of a torus-fixed point on the moduli space
in terms of a choice of a covering by hexagons
of the fundamental region of a real 2-torus
due to Nakamura \cite{Nakamura_HSAGO}.
Many of the arguments are similar
to those in \cite{Ishii_RMMQ}.
There is also a physics paper
by Franco and Vegh \cite{Franco-Vegh_MSGTDM}
which deals with the relation
between brane tilings
and moduli spaces.

{\bf Acknowledgment:}
We thank Alastair King
for a number of very useful remarks,
and the anonymous referee for suggesting several improvements.
This work is supported by Grants-in-Aid for Scientific Research
(No.18540034) and Grant-in-Aid for Young Scientists (No.18840029).

\section{Dimer models and quivers} \label{section:def}

Let $T = \bR^2 / \bZ^2$
be a real two-torus
equipped with an orientation.
A {\em bipartite graph} on $T$
consists of
\begin{itemize}
 \item a set $B \subset T$ of black vertices,
 \item a set $W \subset T$ of white vertices, and
 \item a set $E$ of edges,
       consisting of embedded closed intervals $e$ on $T$
       such that one boundary of $e$ belongs to $B$
       and the other boundary belongs to $W$.
       We assume that two edges intersect
       only at the boundaries.
\end{itemize}
A bipartite graph on $T$ is called a {\em dimer model}
if the set of edges divide $T$ into simply-connected polygons.

A {\em quiver} consists of
\begin{itemize}
 \item a set $V$ of vertices,
 \item a set $A$ of arrows, and
 \item two maps $s, t: A \to V$ from $A$ to $V$.
\end{itemize}
For an arrow $a \in A$,
$s(a)$ and $t(a)$
are said to be the {\em source}
and the {\em target} of $a$
respectively.
A {\em path} on a quiver
is an ordered set of arrows
$(a_n, a_{n-1}, \dots, a_{1})$
such that $s(a_{i+1}) = t(a_i)$
for $i=1, \dots, n-1$.
We also allow for a path of length zero,
starting and ending at the same vertex.
The {\em path algebra} $\bC Q$
of a quiver $Q = (V, A, s, t)$
is the algebra
spanned by the set of paths
as a vector space,
and the multiplication is defined
by the concatenation of paths;
$$
 (b_m, \dots, b_1) \cdot (a_n, \dots, a_1)
  = \begin{cases}
     (b_m, \dots, b_1, a_n, \dots, a_1) & s(b_1) = t(a_n), \\
      0 & \text{otherwise}.
    \end{cases}
$$
A {\em quiver with relations}
is a pair of a quiver
and a two-sided ideal $\scI$
of its path algebra.
For a quiver $\Gamma = (Q, \scI)$
with relations,
its path algebra $\bC \Gamma$ is defined as
the quotient algebra $\bC Q / \scI$.

A dimer model $(B, W, E)$ encodes
the information of a quiver
$\Gamma = (V, A, s, t, \scI)$
with relations
in the following way:
The set $V$ of vertices
is the set of connected components
of the complement
$
 T \setminus (\bigcup_{e \in E} e),
$
and
the set $A$ of arrows
is the set $E$ of edges of the graph.
The directions of the arrows are determined
by the colors of the vertices of the graph,
so that the white vertex $w \in W$ is on the right
of the arrow.
In other words,
the quiver is the dual graph of the dimer model
equipped with an orientation given by
rotating the white-to-black flow on the edges of the dimer model
by minus 90 degrees.

The relations of the quiver are described as follows:
For an arrow $a \in A$,
there exist two paths $p_+(a)$
and $p_-(a)$
from $t(a)$ to $s(a)$,
the former going around the white vertex
connected to $a \in E = A$ clockwise
and the latter going around the black vertex
connected to $a$ counterclockwise.
Then the ideal $\scI$
of the path algebra is
generated by $p_+(a) - p_-(a)$
for all $a \in A$.

A {\em representation} of $\Gamma$ is a module
over the path algebra $\bC  \Gamma$ with relations.
In other words,
a representation is a collection
$((V_v)_{v \in V}, (\psi(a))_{a \in A})$
of vector spaces $V_v$ for $v \in V$
and linear maps $\psi(a) : V_{s(a)} \to V_{t(a)}$
for $a \in A$ satisfying relations in $\scI$.
The Grothendieck group of
the abelian category $\module \bC \Gamma$
of finite dimensional representations of $\bC \Gamma$ is a free abelian group
generated by simple representations
corresponding to the idempotents of $\bC \Gamma$
given as the paths of length zero.
A simple representation corresponding to a vertex $v \in V$
has $V_v = \bC$, $V_w = 0$ for $w \ne v$
and $\psi(a) = 0$ for any $a \in A$.
Let $N$ be the number of vertices of $\Gamma$.
Then the Grothendieck group
is isomorphic to $\bZ^N$
with respect to this basis,
and the class of a module in the Grothendieck group
considered as an element of $\bZ^N$
is called its {\em dimension vector}.
The dimension vector of a representation
$((V_v)_{v \in V}, (\psi(a))_{a \in A})$
is given by $(\dim V_v)_{v \in V}$.

The {\em double} $\overline{Q}$
of a quiver $Q$ is obtained from $Q$
by adding an arrow $\overline{a}$
with $s(\overline{a}) = t(a)$ and $t(\overline{a}) = s(a)$
for each arrow $a$ of $Q$.
A representation $\Psi = ((V_v)_{v \in V}, (\psi(a))_{a \in A})$
of $Q$
such that all $\psi(a)$ are linear isomorphisms
determines a representation $\overline{\Psi}$ of $\overline{Q}$ by
$
 \overline{\Psi}(a) = \Psi(a)
$
and
$
 \overline{\Psi}(\overline{a}) = \Psi(a)^{-1}.
$

A {\em perfect matching}
(or a {\em dimer configuration})
on a dimer model $G = (B, W, E)$
is a subset $D$ of $E$
such that for any vertex $v \in B \cup W$,
there is a unique edge $e \in D$
connected to $v$.
Consider the bipartite graph $\Gtilde$ on $\bR^2$
obtained from $G$ by pulling-back
by the natural projection $\bR^2 \to T$,
and identify the set of perfect matchings of $G$
with the set of periodic perfect matchings of $\Gtilde$.
Fix a reference perfect matching $D_0$.
Then for any
perfect matching $D$,
the union $D \cup D_0$
divides $\bR^2$ into connected components.
The height function $h_{D, D_0}$ is
a locally-constant function on
$\bR^2 \setminus (D \cup D_0)$
which increases (resp. decreases)
by $1$
when one crosses an edge $e \in D$
with the black (resp. white) vertex
on his right
or an edge $e \in D_0$
with the white (resp. black) vertex
on his right.
This rule determines the height function
up to additions of constants.
The height function may not be periodic
even if $D$ and $D_0$ are periodic,
and the {\em height change}
$h(D, D_0) = (h_x(D, D_0), h_y(D, D_0)) \in \bZ^2$
of $D$ with respect to $D_0$
is defined as the difference
\begin{align*}
 h_x(D, D_0) &= h_{D, D_0}(p+(1,0)) - h_{D, D_0}(p), \\
 h_y(D, D_0) &= h_{D, D_0}(p+(0,1)) - h_{D, D_0}(p) 
\end{align*}
of the height function,
which does not depend on the choice of
$p \in \bR^2 \setminus (D \cup D_0)$.
More invariantly,
height changes can be considered
as an element of $H^1(T, \bZ)$.
The dependence of the height change
on the choice of the reference matching
is given by
$$
 h(D, D_1) = h(D, D_0) - h(D_1, D_0)
$$
for any three perfect matchings $D$, $D_0$ and $D_1$.
We will often suppress the dependence of the height difference
on the reference matching
and just write $h(D) = h(D, D_0)$.

For a fixed reference matching $D_0$,
the characteristic polynomial of $G$ is defined by
$$
 Z(x, y)
  = \sum_{D \in \perfmat(G)}
     x^{h_x(D)} y^{h_y(D)},
$$
where $\perfmat(G)$ denotes the set of perfect matchings of $G$.
It is a Laurent polynomial in two variables,
whose Newton polygon coincides
with the convex hull of the set
$$
 \{ (h_x(D), h_y(D)) \in \bZ^2
      \mid D \in \perfmat(G) \}
$$
consisting of height changes of perfect matchings of the dimer model $G$.

A perfect matching can be considered as a set of walls
which block some of the arrows;
for a perfect matching $D$,
let $Q_D$ be the subquiver of $Q$
whose set of vertices is the same as $Q$
and whose set of arrows consists of $A \setminus D$
(recall that $A = E$).
The path algebra $\bC Q_D$ of $Q_D$ is a subalgebra of $\bC Q$,
and the ideal $\scI$ of $\bC Q$
defines an ideal $\scI_D = \scI \cap \bC Q_D$ of $\bC Q_D$.
A path $p \in \bC Q$ is said to be an allowed path
with respect to $D$ if $p \in \bC Q_D$.

%

\begin{figure}[htbp]
\centering
\begin{minipage}{.4 \linewidth}
\centering
\input{conifold_bt.pst}
\caption{A dimer model}
\label{fg:conifold_bt}
\end{minipage}
\begin{minipage}{.45 \linewidth}
\centering
\input{conifold_quiver.pst}
\caption{The corresponding quiver}
\label{fg:conifold_quiver}
\end{minipage}
\end{figure}

\begin{figure}[htbp]
\begin{minipage}{\linewidth}
\centering
\input{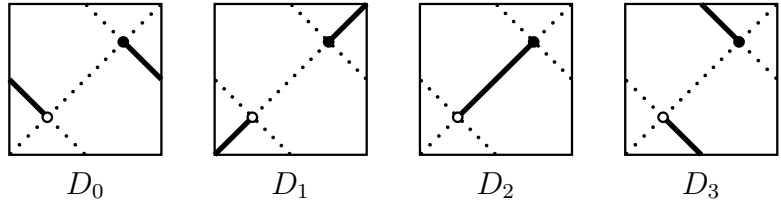}
\caption{Four perfect matchings}
\label{fg:conifold_dimers}
\end{minipage}
\end{figure}

As an example, consider the dimer model
in Figure \ref{fg:conifold_bt}.
The corresponding quiver is shown in Figure \ref{fg:conifold_quiver},
whose relations are given by
$$
 \scI = (d b c - c b d, d a c - c a d, a d b - b d a, a c b - b c a).
$$
This dimer model is non-degenerate,
and has four perfect matchings $D_0, \dots, D_3$
shown in Figure \ref{fg:conifold_dimers}.
Their height changes with respect to $D_0$ are given by
$$
 h(D_0) = (0, 0), \quad
 h(D_1) = (1,0), \quad
 h(D_2) = (0, 1), \quad
 h(D_3) = (1, 1),
$$
so that the characteristic polynomial is
$$
 Z(x, y) = 1 + x + y + x y.
$$

\section{Torus actions on moduli spaces}

Let $\Gamma = (V, A, s, t, \scI)$ be the quiver with relations
obtained from a dimer model $(B, W, E)$
on a real two-torus $T = \bR^2 / \bZ^2$ and
$\scMtilde$ be the set of
representations of $\Gamma$ with dimension vector $(1, \dots, 1)$.
In other words,
$\scMtilde$ is the subset of $\bC^A$
consisting of linear maps
$\psi(a) : \bC s(a) \to \bC t(a)$
for arrows $a \in A$
satisfying relations.
$\scMtilde$ has a natural scheme structure as a closed subscheme
of $\bC^A$ defined by the ideal generated by the relations.
Let $\prod_v \Aut(V_v) \cong (\bCx)^V$
be the product of the set of automorphisms
of the vector spaces attached to vertices
of the quiver.
Two representations $((V_v)_{v \in V}, (\psi(a))_{a \in A})$
and $((W_v)_{v \in V}, (\phi(a))_{a \in A})$
will be called isomorphic
if there is an element $(g_v)_{v \in V}$
such that for any $a \in A$,
the following diagram commutes;
$$
\begin{CD}
 V_{s(a)} @>{\psi(a)}>> V_{t(a)} \\
 @V{g_{s(a)}}VV @VV{g_{t(a)}}V \\
 W_{s(a)} @>{\phi(a)}>> W_{t(a)}.
\end{CD}
$$
The diagonal subgroup $\bCx \subset (\bCx)^V$ acts trivially on $\scMtilde$
and the quotient $\scG=(\bCx)^V/\bCx$ acts faithfully on $\scMtilde$.
The set-theoretic quotient
of $\scMtilde$ by the action of $\scG$
will be denoted by $\scM$.
Let $\bT \subset \scM$ be the subset
consisting of isomorphism classes
$[(\psi(a))_{a \in A}]$
such that $\psi(a) \in \bCx$ for any $a \in A$.
It has a structure of an algebraic torus
by the pointwise multiplication:
For two elements
$[(\psi(a))_{a \in A}]$ and $[(\phi(a))_{a \in A}]$ of $\bT$,
their composition is defined by
\begin{equation} \label{eq:pointwise_multiplication}
 [(\psi(a))_{a \in A}] \cdot [(\phi(a))_{a \in A}]
   = [(\psi(a) \cdot \phi(a))_{a \in A}], 
\end{equation}
which gives an element of $\bT$ again.

The set $\scM$ of isomorphism classes
does not have a good geometric structure.
We use the notion of stability
introduced by King \cite{King}
in order to construct moduli schemes of quiver representations.

\begin{dfn}
For $\theta \in \Hom_\bZ(K(\module \bC \Gamma), \bZ)$
such that $\theta((1, \dots, 1)) = 0$,
a representation $\Psi \in \scMtilde$ is said to be {\em $\theta$-stable}
if for any non-trivial subrepresentation
(i.e., subobject in $\module \bC \Gamma$)
$S \subsetneq \Psi$,
we have $\theta(S) >0$.
$\Psi$ is {\em $\theta$-semistable}
if $\theta(S) \ge 0$ holds instead of $\theta(S)>0$ above.
\end{dfn}

In this definition, $\theta$ corresponds to a character of $\scG$.
King \cite{King} proved that
there is a moduli scheme $\overline{\scM_\theta}$
which parameterizes the S-equivalence classes
of $\theta$-semistable representations in $\scMtilde$.
It contains the moduli scheme $\scM_\theta$
which parameterizes the isomorphism classes of
$\theta$-stable representations as an open set.

Note that if $\psi(a) \in \bCx$ for any $a \in A$,
then $\Psi$ doesn't have any non-trivial subrepresentation
and hence is $\theta$-stable
for any $\theta \in \Hom_\bZ(K(\module \bC \Gamma), \bZ)$.
Thus $\bT$ is naturally contained in $\scM_\theta$ for any $\theta$.
Moreover, there is an action of $\bT$ on $\scM_\theta$
defined by the pointwise multiplication
just as in (\ref{eq:pointwise_multiplication}).

Now consider the complex
$$
 1 \to (\bCx)^V \xto{d^1} (\bCx)^A \xto{d^2} (\bCx)^F \to 1.
$$
Here, $F = B \cup W$ is the set of vertices
of the dimer model
which is in one-to-one correspondence
with the set of faces of the quiver.
The map $d^1$ is defined by
$$
\begin{array}{rccc}
 d^1 : & (\bCx)^V & \rightarrow & (\bCx)^A \\
   & \rotatebox{90}{$\in$} & & \rotatebox{90}{$\in$} \\
   & (g_v)_{v \in V} & \mapsto
              & (g_{s(a)}^{-1} \cdot g_{t(a)} )_{a \in A},
\end{array}
$$
and the map $d^2$ sends $(\psi(a))_{a \in A}$
to $(\phi_f)_{f \in F}$,
where $\phi_f$ is the product of all $\psi(a)$
such that $a \in A = E$ is connected to $f \in B \cup W$.
The above complex is the cochain complex
computing the $\bCx$-valued cohomologies of $T$
with respect to the polygonal division of $T$
determined by the quiver $\Gamma$.

Let $\bT_0$ denote the diagonal subgroup of $(\bCx)^F$
and $\bTtilde \subset \scMtilde$ the preimage of $\bT \subset \scM$.
Then one has
$$
 \bTtilde = (d^2)^{-1}(\bT_0),
$$
and the cohomology group in the middle
of the subcomplex
$$
 1 \to (\bCx)^V \xto{d^1} \bTtilde \xto{d^2} \bT_0 \to 1
$$
is isomorphic to $H^1(T, \bCx)$.
It follows from the definition that
$$
 \bT = \bTtilde / \Image d^1,
$$
and hence one has an exact sequence
\begin{equation} \label{eq:bT_exact}
 1 \to H^1(T, \bCx) \to \bT \to \bT_0. 
\end{equation}
This proves the following:
\begin{lmm}
The dimension of the algebraic torus $\bT$ is either two or three. 
\end{lmm}

\section{Coordinates around $\bT$-fixed points}\label{section:coordinates}

Suppose that a representation
$
 \Psi = (\psi(a))_{a \in A} \in \scMtilde
$
represents a point $[\Psi] \in \scM_{\theta}$,
which is fixed by the action of $\bT$.
Let $\Gamma_\Psi$ be the subquiver of $\Gamma$
whose set of vertices is $V$ and
whose set of arrows consists of arrows $a \in A$
such that $\psi(a) \ne 0$.
The stability of $\Psi$ implies that $\Gamma_\Psi$ is connected.
Moreover, we have the following.
\begin{lmm}\label{lemma:lift}
If $[\Psi]\in \scM$ is fixed by the action of $H^1(T, \bCx) \subset \bT$,
then $\Gamma_\Psi$ can be lifted to a subquiver $\Gamma'_\Psi$ of $\Gammatilde$, which is
isomorphically mapped to $\Gamma_\Psi$.
$\Gamma'_\Psi$ is unique up to translations by $\bZ^2 \subset \bR^2$.
\end{lmm}
\begin{proof}
Fix a vertex $v_0$ of $\Gamma_\Psi$ and lift it to a vertex $\vtilde_0$ of $\Gammatilde$.
For a vertex $v$ of $\Gamma_\Psi$, take a path $p$ of the double $\overline{\Gamma_\Psi}$ of
$\Gamma_\Psi$ starting from $v_0$ and ending at $v$.
We can lift $p$ to a path $\ptilde$ of the double of ${\Gammatilde}$ starting from $\vtilde_0$.
We will show that the end point of $\ptilde$ does not depend on the choice of $p$.

Assume that there are two paths $p_1$ and $p_2$ of $\overline{\Gamma_\Psi}$ starting from $v_0$ and ending at $v$
such that the endpoints of their lifts $\ptilde_1$ and $\ptilde_2$ are different.
The path $\gamma:=p_2 \cdot (p_1)^{-1}$ is a loop starting from $v$ and the assumption
implies that it determines a non-trivial class $[\gamma] \in H_1(T, \bZ)$.
Consider the value $\psi(\gamma)$ of $\Psi$ at $\gamma$;
we can define values of $\Psi$ for arrows and paths of the double $\overline{\Gamma_\Psi}$
in an obvious way.
Since $[\gamma]$ is a non-trivial class, there is $g \in H^1(T, \bCx)$ with
$(g\cdot \psi)(\gamma) \ne \psi(\gamma)$.
This contradicts the assumption that $H^1(T, \bCx)$ fixes $[\Psi]$.
\end{proof}


Let $F_\Psi$ be the closure of the union of the connected components of
$
 \bR^2 \setminus (\bigcup_{e \in \Etilde} e)
$
corresponding to the vertices of $\Gamma'_\Psi$.
It is a fundamental domain for the action of $\bZ^2$ on $\bR^2$.

Now recall that the set $A$ of arrows of the quiver
is identified with the set $E$ of the edges of the dimer model.
Thus for an edge $e \in E$,
we write $\psi(e)$ for the value of $\psi$ at the arrow corresponding to $e$.
For an edge $e \in \Etilde$,
we also write $\psi(e)$
for the value of $\psi$ at the corresponding edge of $E$.
Let $\delta_\Psi$ be the union of edges $e \in \Etilde$
satisfying the following:
\begin{itemize}
\item $\psi(e)=0$
\item $e$ intersects with another edge $e'$ with $\psi(e')=0$.
\end{itemize}
The boundary of the fundamental domain $F_\Psi$ is obviously contained in $\delta_\Psi$.
On the other hand,
the relations of the quiver imply that
there are no end points in $\delta_\Psi$,
and therefore the interior of $F_\Psi$
does not intersect with $\delta_\Psi$.
Thus $\delta_\Psi$ is the union of
the translations of the boundary $\partial F_\Psi$ of $F_\Psi$:
\begin{equation}\label{equation:boundary}
\delta_\Psi = \bigcup_{m \in \bZ^2} (\partial F_\Psi + m)
\end{equation}

\begin{lmm}\label{lemma:=1}
By replacing $\Psi$ with a representation equivalent to $\Psi$,
we may assume $\psi(a)=1$ for all arrows $a$ of $\Gamma_\Psi$.
\end{lmm}
\begin{proof}
The assertion means that
we can attach a complex number $g_v \in \bCx$
to each vertex $v \in V$ such that
$\psi(a)=g_{t(a)}g_{s(a)}^{-1}$ for any arrow $a$ in $\Gamma_\Psi$.
Fix a vertex $v_0 \in V$.
For any $v \in V$,
we take a path $p$ in the double quiver $\overline{\Gamma_\Psi}$
starting from $v_0$ and ending at $v$ and we want to put $g_v = \psi(p)$.
If we show that $\psi(p)$ does not depend on the choice of $p$, we are done.
Take two such paths $p_1$ and $p_2$.
Lemma \ref{lemma:lift} implies that $p_1$ and $p_2$ are homotopic in $T$
and \eqref{equation:boundary} shows that the homotopy is generated by the relations
\begin{itemize}
\item
$p_+(a) \sim p_-(a)$ for arrows $a$ of $\Gamma$ such that $p_+(a)$ is a path in $\Gamma_\Psi$.
\item
$a^{-1}\cdot a \sim e_{s(a)}$ and $a \cdot a^{-1} \sim e_{t(a)}$ for arrows $a$ of ${\Gamma_\Psi}$.
\end{itemize}
Thus we obtain $\psi(p_1)=\psi(p_2)$.
\end{proof}

From now on,
we assume $\psi(a)=1$ for all arrows $a$ of $\Gamma_\Psi$.
In other words, $\psi(a)$ is either $0$ or $1$.
Consider the following subset $U_\Psi$ of $\scMtilde$:
$$
 U_\Psi =
   \{ \Phi = (\phi(a))_{a \in A} \in \scMtilde
       \mid \phi(a) = 1 \text{ if } \psi(a) = 1 \}.
$$
$U_\Psi$ is naturally a closed subscheme of $\scMtilde$.

\begin{lmm}\label{lemma:open_immersion}
Every point in $U_\Psi$ is $\theta$-stable and the natural morphism
$U_\Psi \to \scM_\theta$ is an open immersion.
Thus $U_\Psi$ can be regarded as a $\bT$-invariant affine open neighborhood of $[\Psi]$
in $\scM_\theta$.
\end{lmm}
\begin{proof}
Suppose $\Phi \in U_\Psi$.
Since the dimension vector of $\Phi$ is $(1,1,\dots,1)$,
a subrepresentation of $\Phi$ is determined by a subset $V'$ of $V$.
By the definition of $U_\Psi$, $V'$ also determines a subrepresentation
of $\Psi$.
Thus the $\theta$-stability for $\Psi$ implies that for $\Phi$; the first assertion follows.
For the second assertion, put
$$
\Utilde_\Psi=\{ \Phi= (\phi(a))_{a \in A}\in \scMtilde \mid \phi(a) \ne 0 \text{ if } \psi(a)=1\}.
$$
This is an open subscheme of $\scMtilde$.
Then the same argument as in Lemma \ref{lemma:=1} shows that the morphism
$$
\scG \times  U_\Psi \to \Utilde_\Psi
$$
induced by the action of $\scG$ on $\Utilde_\Psi$ is an isomorphism.
Thus $U_\Psi$ is a section of the morphism from $\Utilde_\Psi$ to its quotient by $\scG$.
\end{proof}

\begin{lmm}\label{lemma:4or6}
Either of the following two cases must occur:
\begin{enumerate}
\item There are four quadrivalent points of the graph $\delta_\Psi$ lying in
$\partial F_\Psi$, and
there are no points of valency three or greater than four.
\item
There are six trivalent points of $\delta_\Psi$ lying in
$\partial F_\Psi$,
and there are no points of valency greater than three.
\end{enumerate}
\end{lmm}
\begin{proof}
Let $a_n$ be the number of points of valency $n$ of $\delta_\Psi$ lying in
$\partial F_\Psi$.
These points divide $\partial F_\Psi$ into
$\left(\sum_{n \ge 3} a_n\right)$ parts so that we can regard $\partial F_\Psi$
as a polygon with $\left(\sum_{n \ge 3} a_n\right)$ edges.
Since a point of valency $n$ is contained in $n$ translations of $F_\Psi$, the equation that the topological Euler number of $T$ is zero leads to
$$
1 - \frac{1}{2}\sum_{n\ge 3} a_n + \sum_{n\ge 3} \frac{a_n}{n} =0.
$$
It is easy to see from this that there are only two possibilities as stated.
\end{proof}

\begin{lmm}\label{lemma:local_chart}
If $\dim \bT=3$, then it holds that $\bT \subset U_\Psi \cong \bC^3$.
If $\dim \bT=2$, then $U_\Psi$ is the disjoint union of $\bT$ and the isolated point $\{[\Psi]\}$.
\end{lmm}
\begin{proof}
We first consider the case 1 of Lemma \ref{lemma:4or6}.
Assume that $v_1, v_2, v_3, v_4$ are the quadrivalent points of $\delta_\Psi$
lying on $\partial F_\Psi$, labeled counterclockwise.
These points are mapped to a common vertex $v \in T$ of the dimer model.
Since $v_1$ is a quadrivalent point of $\delta_\Psi$,
there are four edges $e_1, e_2, e_3, e_4$ of $\Etilde$
that are connected to $v_1$ and that satisfy $\psi(e_i)=0$.
The four points $v_i$ divides $\partial F_\Psi$ into four parts and we may assume that $e_1$ is
on the part between $v_1$ and $v_2$,
and $e_2$ is between $v_1$ and $v_4$.
We further assume $e_1, e_2, e_3, e_4$ are arranged counterclockwise
around $v_1$.

Now take $\Phi \in U_\Psi$ and put $t_i=\phi(e_i)$.
Then, the relations $p_+(a)=p_-(a)$ for arrows $a$ corresponding to the edges in
$\partial F_\Psi$ determine the values of $\Phi$ at the edges on $\partial F_\Psi$:
\begin{itemize}
\item
For an edge $e$ of $\partial F_\Psi$
between $v_1$ and $v_2$,
$\phi(e)$ coincides with either $t_1$ or $t_2t_3t_4$,
depending on the configuration of the colors of the vertices of $e$.
In particular,
since the colors of $v_1$ and $v_2$ are the same,
we have $t_3=t_2t_3t_4$ and $t_1=t_4t_1t_2$.
\item
Similarly, for an edge $e$ of $\partial F_\Psi$
between $v_1$ and $v_4$, $\phi(e)$ coincides with either $t_2$ or $t_3t_4t_1$,
and we have $t_4=t_3t_4t_1$ and $t_2=t_1t_2t_3$.
\end{itemize}
Moreover, since $a p_+(a)$ does not depend on an arrow $a$, we obtain the following:
\begin{itemize}
\item For an edge $e \in \Etilde$ with $\psi(e) =0$ that is not in $\delta_\Psi$,
we must have $\phi(e)=t_1t_2t_3t_4$.
\end{itemize}
Thus $\Phi \in U_\Psi$ is determined by the point $(t_1, t_2, t_3, t_4) \in \bC^4$.
Conversely, for any point in $\bC^4$ that satisfies the relations $t_3=t_2t_3t_4$, $t_1=t_4t_1t_2$, $t_4=t_3t_4t_1$ and $t_2=t_1t_2t_3$, we can find a corresponding point in $U_\Psi$.
Solving these four equations, we obtain
$$
U_\Psi \cong \{(t_1, t_2, t_3, t_4) \in \bC^4 \mid (t_1, t_2, t_3, t_4) =0 \text{ or } t_1t_3 = t_2 t_4 =1 \}.
$$
The two-dimensional component defined by $t_1t_3=t_2t_4=1$ is $\bT$-invariant and is contained in $\bT$;
hence it coincides with $\bT$ which must be two-dimensional.
The origin of $\bC^4$ corresponds to $[\Psi]$.

Next we consider the case 2 of Lemma \ref{lemma:4or6}.
Let $v_1, \dots, v_6$ denote the six trivalent points of $\delta_\Psi$
lying counterclockwise on $\partial F_\Psi$.
In this case, $v_1, v_3$ and  $v_5$ are in a single $\bZ^2$-orbit in $\bR^2$, and  $v_2, v_4$ and $v_6$ are in another orbit.
$v_1$ is connected to three edges $e_1, e_2, e_3$ of $\Etilde$ that satisfy $\psi(e_i)=0$
and $v_2$ is connected to $e_4, e_5,e_6$ similarly.
We may assume that $e_1$ and $e_4$ are on the part of $\partial F_\Psi$ cut out by $v_1$ and $v_2$
which contains no other $v_i$,
and that $e_1, e_2, e_3$ and $e_4, e_5, e_6$ are arranged counterclockwise
around $v_1$ and $v_2$ respectively.
As in the case 1, $\Psi$ is determined by $t_i:=\phi(e_i)$($i=1,\dots,6$)
and we can see
\begin{itemize}
\item If one of $v_1$ and $v_2$ is black and the other is white,
then $(t_1, t_2, t_3)=(t_4, t_5, t_6)$.
\item If the colors of $v_1$ and $v_2$ are the same,
then $t_i=t_jt_k$ for $(i,j,k)=(1,5,6)$, $(2,6,4)$, $(3,4,5)$, $(4,2,3)$, $(5,3,1)$ and $(6,1,2)$.
\end{itemize}
In the first case, $(t_1, t_2, t_3)$ gives rise to an isomorphism $U_\Psi \cong \bC^3$ and
$\bT$ coincides with the open subset defined by $t_1t_2t_3 \ne 0$.
In the second case, we can see
$$
U_\Psi \cong \{(t_1, t_2, t_3) \in \bC^3 \mid (t_1, t_2, t_3)=0 \text{ or } t_1t_2t_3 =1 \},
$$
which, as in the case 1, is the union of $\{[\Psi]\}$ and the two-dimensional tours $\bT$.
\end{proof}

\section{Moduli spaces as crepant resolutions} \label{section:resolution}

Our definition of dimer models
in \S \ref{section:def} contains
a lot of ``inconsistent'' ones
from a physics point of view
(see Hanany--Vegh \cite{Hanany-Vegh_QTBR}).
Here we introduce a condition which should be
necessary (but not sufficient) for the consistency.

A dimer model is said to be {\em non-degenerate}
if for any edge $e \in E$,
there exists a perfect matching $D$
such that $e \in D$.
An {\em R-charge} on a dimer model $G = (B, W, E)$
is a collection of positive real numbers
$R_e \in \bR_{>0}$ indexed by edges $e \in E$,
satisfying
\begin{equation} \label{R-charge}
 \sum_{\substack{e \in E \\ e \ni v }} R_e = 2
\end{equation}
for each vertex $v \in B \cup W$.
If $G$ is non-degenerate,
one can define an R-charge
by averaging
$$
 R_e = \frac{2}{|\perfmat(G)|} \sum_{D \in \perfmat(G)} \chi_D(e)
$$
over the set $\perfmat(G)$ of perfect matchings.
Here $\chi_D$ is the characteristic function
of the subset $D \subset E$.

\begin{rem}
Alastair King pointed to us that
the Birkhoff--von Neumann theorem implies that
the non-degeneracy condition is in fact equivalent
to the existence of an R-charge.
He also remarked that
Hall's marriage theorem implies that
this condition is also equivalent to the following
{\em strong marriage condition};
every proper subset of blacks is connected to strictly more whites
and vice versa.
\end{rem}


Take the parameter
$0 \in \Hom_\bZ(K(\module \bC \Gamma), \bZ)$
and consider the corresponding moduli space $\overline{\scM_0}$.
Since any representation of $\Gamma$ is $0$-semistable, this is the categorical
quotient of $\scMtilde$ by the action of $\scG$.
Hence $\overline{\scM_0}$ is an affine scheme with a distinguished point
$[0] \in \overline{\scM_0}$ which is the image of $0 \in \scMtilde \subset \bC^A$.
Moreover, for any parameter $\theta$, we have a projective morphism
$\overline{\scM_\theta} \to \overline{\scM_0}$.

\begin{prp} \label{prop:crepant}
Let $(B, W, E)$ be a non-degenerate dimer model.
Then we have $\dim \bT=3$, and
for a generic parameter $\theta$,
the moduli space $\scM_\theta$
is smooth and irreducible with the trivial canonical bundle $K_{\scM_\theta}$.
\end{prp}

\begin{proof}
We may assume the existence of an R-charge
$(R_e)_e \in (\bQ_{>0})^E$
satisfying \eqref{R-charge}.
Take a positive integer $N$ such that $r_e:=N R_e$ is an integer.
Then $t \mapsto (t^{r_e})_{e \in E}$ is a one parameter subgroup of $\bT$
not contained in $H^1(T, \bCx)$,
and hence we have $\dim \bT=3$.

Take an arbitrary point $[\Phi]=(\phi(a)_{a\in A})$ in $\scM_\theta$.
We will show that there is a $\bT$-fixed point $[\Psi]\in \scM_\theta$ such that
$[\Phi]\in U_\Psi$.

Consider the morphism $\xi: \Spec \bC[t] \to \overline{\scM_0}$
defined by $t \mapsto [(t^{r_a}\phi(a))_{a\in A}]$.
We have $\xi(1)=[\Phi]$ and $\xi(0)=[0]$.
Moreover, for $t \ne 0$, $\xi(t)$ is $\theta$-stable.
By virtue of the valuative criterion
for the projective morphism $\scM_\theta \to \overline{\scM_0}$,
we can lift $\xi$ to $\widetilde{\xi}: \Spec \bC[t] \to \scM_\theta$.
Since $U_\Psi$ is a $\bT$-invariant open subset and
$\widetilde{\xi}(\Spec \bC[t] \setminus \{0\})$ is contained
in a single $\bT$-orbit,
it suffices to show $\widetilde{\xi}(0) \in U_\Psi$
for some $\bT$-fixed point $[\Psi]$.
Since the fiber of
$
 \scM_\theta \to \overline{\scM_0}$
over
$
 [0] \in \overline{\scM_0}
$
is a $\bT$-invariant closed subscheme
projective over $\Spec \bC$,
we can find such $\Psi$
as a limit point of the $\bT$-action.
Hence we have $[\Phi] \in U_\Psi$, where $U_\Psi$ contains $\bT$ and is isomorphic to $\bC^3$ by Lemma \ref{lemma:local_chart}.
Since $[\Phi]$ is arbitrary, $\scM_\theta$ is smooth and irreducible.

Now we prove that the canonical bundle of $\scM_\theta$ is trivial.
As in the proof of Lemma \ref{lemma:local_chart},
we have a coordinate $(t_1, t_2, t_3)$ on $U_\Psi$.
We show that we can patch the $3$-forms
$dt_1 \wedge dt_2 \wedge dt_3$ on $U_\Psi$
to obtain a global $3$-form on $\scM_\theta$.
Let $[\Psi]$ and $[\Phi]$ be two $\bT$-fixed points
on $\scM_\theta$.
Then we have coordinates $t_1, t_2, t_3$ on $U_\Psi$
and $s_1, s_2, s_3$ on $U_\Phi$ respectively.
On the torus $\bT$, we can express $t_1, t_2, t_3$
as Laurent monomials in $s_1, s_2, s_3$, and vice versa.
Thus $t_1, t_2, t_3$ and $s_1, s_2, s_3$ are related
by a matrix in $\GL(3, \bZ)$.
Moreover, we have $t_1 t_2 t_3 = s_1 s_2 s_3$.
These two facts imply
$
 dt_1\wedge dt_2 \wedge dt_3 = \pm ds_1 \wedge ds_2 \wedge ds_3.
$
To determine the sign,
recall that the edges $e_1, e_2, e_3$
that correspond to $t_1, t_2, t_3$
in the proof of Lemma \ref{lemma:local_chart}
are arranged counterclockwise.
We can make the same assumption on the choice of $s_1, s_2, s_3$.
Then we can see that
the matrix in $\GL(3, \bZ)$ has the determinant one,
so that $dt_1\wedge dt_2 \wedge dt_3 = ds_1 \wedge ds_2 \wedge ds_3$.
\end{proof}

\section{Perfect matchings and toric divisors on moduli spaces}

In this section, we discuss the relation between perfect matchings
and $\bT$-invariant divisors on moduli spaces.
Throughout this section,
we assume that
$G=(B, W, E)$ is a non-degenerate dimer model.

For a generic $\theta$
and a two-dimensional $\bT$-orbit $Z$ in $\scM_\theta$,
pick a representation $\Psi = [(\psi(a))_{a \in A}] \in Z$ and put
$$
D_Z=\{a \in A \mid \psi(a)=0 \}. 
$$
This does not depend on the choice of $\Psi$ in $Z$.

\begin{lmm}\label{lemma:2-dim_orbit_1}
If $\theta$ is generic, then $D_Z$ is a perfect matching for any two-dimensional $\bT$-orbit $Z$
in $\scM_\theta$.
\end{lmm}
\begin{proof}
Take a $\bT$-fixed point $[\Phi] \in \overline{Z}$ and
consider the affine open neighborhood $U_\Phi$ of $[\Phi]$
appearing in \S \ref{section:coordinates}.
As in the proof of Lemma \ref{lemma:local_chart}, there is a coordinate $(t_1, t_2, t_3)$ on $U_\Phi$ that gives rise to an isomorphism $U_\Phi \cong \bC^3$.
The action of $\bT$ on $U_\Psi$ is diagonalized
with respect to this coordinate and
hence $Z \subset U_\Phi$ is defined by $t_i=0$ and
$t_j t_k \ne 0$ where $\{i, j, k\}=\{1, 2, 3\}$.
Then it follows from the proof of Lemma \ref{lemma:local_chart}
that $D_Z$ is a perfect matching.
\end{proof}

Suppose that $D$ is a perfect matching.
For $t \in \bC$, we define $\Psi_t=(\psi_t(a))_{a \in A}$ by
\begin{equation}\label{equation:split}
\psi_t(a)= \begin{cases} t &\text{if } a \in D \subset E=A \\
1 &\text{if } a \nin D \end{cases}.
\end{equation}
Then we can see that $\Psi_t$ satisfies the relation of the quiver
and the graph $\Gamma_{\Psi_t}$ is connected.

\begin{lmm}\label{lemma:2-dim_orbit_2}
There is a generic parameter $\theta$  such that $\Psi_0$ is $\theta$-stable.
Moreover, the $\bT$-orbit $Z$ of $[\Psi_0]$ in $\scM_\theta$ is two-dimensional and it satisfies
$D=D_Z$.
\end{lmm}
\begin{proof}
To find a parameter $\theta$ such that $\Psi_0$ is $\theta$-stable,
we can use an idea from Sardo--Infirri \cite{Sardo-Infirri}:
For an arrow $a \in A \setminus D$,
take an arbitrary positive rational number $\xi_a$.
For a vertex $v$ of the quiver, we put
$$
\theta(v)= \sum_{\substack{a \in A \setminus D \\ t(a)=v}} \xi_a - \sum_{\substack{a \in A \setminus D \\ s(a)=v}} \xi_a.
$$
Then, for any non-trivial subrepresentation $S$ of $\Psi_t$, we have
$$
\theta(S)= \sum_{\substack{a \in A \setminus D \\ s(a) \nin S \\ t(a) \in S}} \xi_a > 0,
$$
which shows that $\Psi_t$ is $\theta$-stable.

For the genericity of $\theta$,
it suffices to show that we can take $\theta$ so that
$$
\theta(S)= \sum_{\substack{a \in A \setminus D \\ s(a) \nin S \\ t(a) \in S}} \xi_a
- \sum_{\substack{a \in A \setminus D \\ s(a) \in S \\ t(a) \nin S}} \xi_a
\ne 0
$$
for an arbitrary non-empty subset $S \subsetneq V$.
This is achieved if $(\xi_a)_{a \in A \setminus D}$ is sufficiently general.

$H^1(T, \bCx)$ acts freely on $Z$
since any element in $H_1(T, \bZ)$ can be represented
by a linear combination of paths of the double $\overline{Q_D}$.
This shows that $Z$ is two-dimensional.
It follows from the definition of $\Psi_0$ that $D=D_Z$.
\end{proof}

Consider the closure $X'$ of $\bT$
in the moduli space $\overline{\scM_0}$
corresponding to the parameter $0$.
Since it is not a priori clear if $X'$ is normal,
we take the normalization $X$ of $X'$
which is an affine toric variety.
Proposition \ref{prop:crepant} is saying that $\scM_\theta$ is a crepant resolution of $X$ for a generic $\theta$.

Now let $\Delta \subset H^1(T, \bZ)$ be
the Newton polygon of the characteristic polynomial
(i.e., the convex hull of height changes) of the dimer model
with respect to any fixed perfect matching.
Then we have the following:

\begin{prp} \label{prop:affine_coordinate_ring}
The affine coordinate ring of $X$ is isomorphic to the semigroup ring $\bC[(\Cone(\Delta \times \{ 1 \}))^{\circ}]$ of the dual cone of the cone over
$\Delta \times \{ 1 \} \subset H^1(T, \bZ) \times \bZ$.
\end{prp}

\begin{proof}
Put $N=\Hom(\bCx, \bT)$.
Then the  affine toric variety $X$ is determined by a cone $C \subset N$.
For a perfect matching $D$ of $G=(B, W, E)$,
let
$$
\iota_D : \bCx \to \bT
$$
be the homomorphism
sending $t \in \bCx$ to the representation
$\Psi_t$ defined by \eqref{equation:split}.
Since $\iota_D(\bCx)$ is the stabilizer of $[\Psi_0] \in Z \subset \scM_\theta$
where the orbit $Z$ of $[\Psi_0]$ is two-dimensional by Lemma \ref{lemma:2-dim_orbit_2},
Lemmas \ref{lemma:2-dim_orbit_1} and \ref{lemma:2-dim_orbit_2} imply that
the cone $C$ is generated by the set
$$
 \{ \iota_D \mid D \in \perfmat(G) \} \subset N.
$$
Recall that we have an exact sequence
$$
 1 \to H^1(T, \bCx) \to \bT \to \bCx \to 1. 
$$
The map $\iota_D$ gives a splitting
$$
 \bT \cong H^1(T, \bCx) \times \iota_D(\bCx)
$$
of the above exact sequence.
Let
$$
 \pi_D : \bT \to H^1(T, \bCx)
$$
be the projection with respect to this splitting.
Now fix a reference perfect matching $D_0$,
and hence the splitting
$
 \bT \cong H^1(T, \bCx) \times \bCx
$
given by $\iota_{D_0}$.
Then under the corresponding splitting
$N \cong H^1(T, \bZ) \times \bZ$, we have
$\iota_D \in  H^1(T, \bZ) \times \{1\} $
for any perfect matching $D$.
Therefore it suffices to show that $\Delta$ is the convex hull of
$$
 \{ \pi_{D_0} \circ \iota_D \mid
     D \in \perfmat(G) \} \subset H^1(T, \bZ),
$$
where we have identified $\Hom(\bCx, H^1(T, \bCx))$ with $H^1(T, \bZ)$.

The projection $\pi_D$ is defined as follows:
For a homology class $C \in H_1(T, \bZ)$,
choose an allowed path $p_C$
with respect to $D$
(i.e., a path which does not contain
any arrow $a \in D \subset E = A$)
whose homology class lies in $C$.
Then for $\psi \in \bT$,
$$
 \pi_{D}(\psi)(C) = \psi(p_C) \in \bCx.
$$
It follows that
the height change $h_{D, D_0} \in H^1(T, \bZ)$
of $D$ with respect to the reference perfect matching $D_0$
coincides with $\pi_{D_0} \circ \iota_D$.
Since $\Delta$ is the convex hull of the set of height changes,
we are done.
\end{proof}

By combining Proposition \ref{prop:crepant}
with Proposition \ref{prop:affine_coordinate_ring},
we obtain the main theorem in this paper:

\begin{thm} \label{th:main}
Let $(B, W, E)$ be a non-degenerate dimer model.
Then for a generic parameter $\theta$,
$\scM_\theta$ is a crepant resolution
of $\Spec \bC[(\Cone(\Delta \times \{ 1 \}))^{\circ}]$.
\end{thm}

For example,
the moduli space $\scM_{\theta}$ for the dimer model
in Figure \ref{fg:conifold_bt}
and a generic stability parameter $\theta$ is the total space
of the direct sum $\scO_{\bP^1}(-1) \oplus \scO_{\bP^1}(-1)$
of the tautological bundle $\scO_{\bP^1}(-1)$ on the projective line $\bP^1$.
There is a real-codimension one wall in the space of stability parameters
and the moduli space flops as one moves from one chamber to the other.


\ \\

\noindent
Akira Ishii

Department of Mathematics,
Graduate School of Science,
Hiroshima University,
1-3-1 Kagamiyama,
Higashi-Hiroshima,
739-8526,
Japan

{\em e-mail address}\ : \ akira@math.sci.hiroshima-u.ac.jp

\ \\

\noindent
Kazushi Ueda

Department of Mathematics,
Graduate School of Science,
Osaka University,
Machikaneyama 1-1,
Toyonaka,
Osaka,
560-0043,
Japan.

{\em e-mail address}\ : \  kazushi@math.sci.osaka-u.ac.jp

\end{document}